\documentclass[a4paper, 12pt]{article}

\usepackage[latin1]{inputenc}
\usepackage[dvips]{graphics}
\usepackage{epsfig}
\usepackage{bbm}
\usepackage{wrapfig}
\usepackage{float}
\usepackage{setspace}
\usepackage{picins}
\usepackage{amssymb}
\usepackage{graphics}
\usepackage{mathrsfs}
\usepackage{amsthm}
\usepackage{amsmath}

\normalsize
\theoremstyle{definition}

\newtheorem{theorem}{Theorem}

\newtheorem{lemma}{Lemma}
\newtheorem{conjecture}{Conjecture}

\newtheorem*{corollary}{Corollary}

\theoremstyle{remark}

\newtheorem*{proofthm1}{Proof of Theorem~\ref{thm:rootedtreealt}}

\title{A proof of the rooted tree alternative conjecture}
\author{Mykhaylo Tyomkyn}

\begin{document}

\maketitle

\begin{abstract}
In~\cite{bonato} Bonato and Tardif conjectured that the number of isomorphism classes of trees mutually embeddable with a given tree $T$ is either $1$ or $\infty$. We prove the analogue of their conjecture for rooted trees. We also discuss the original conjecture for locally finite trees and state some new conjectures.
\end{abstract}

\section{Introduction}
\label{sec:introduction}

Embeddings and isomorphisms of infinite graphs have been much studied since Rado's~\cite{rado} classical paper in~1964.
We say that a graph $G$ \emph{embeds} into a graph $H$ if there exists an injective graph homomorphism $\phi \colon G \rightarrow H$. Equivalently, $G$ embeds into $H$ if $G$ is isomorphic to a subgraph of $H$. We say that $G$ \emph{strongly embeds} into $H$ if $G$ is isomorphic to an \emph{induced} subgraph of $H$. Let us call two graphs $G$ and $H$ \emph{twins} or \emph{mutually embeddable}, in notation $G \sim H$, if $G$ embeds into $H$ and vice versa. Similarly, we call $G$ and $H$ \emph{mutually strongly embeddable}, written $G \approx H$, if $G$ strongly embeds into $H$ and vice versa. It is clear that $\sim$ and $\approx$ are equivalence relations and that $\sim$ is refined by $\approx$ which, in turn, is refined by the isomorphism relation $\cong$.

If $G$ and $H$ are two finite mutually embeddable graphs, then they are isomorphic. However, for infinite graphs this is no longer the case. For example, if a graph $G$ contains a copy of the complete countably infinite graph $K_{\omega}$ then $G$ and $K_{\omega}$ are twins. Also, if a countable graph $G$ contains the Rado graph $R$ (see \cite{diestel}, \cite{rado}) as an induced subgraph, then $G$ and $R$ are mutually strongly embeddable, since $R$ contains \emph{every} countable graph as an induced subgraph. In both examples the number of isomorphism classes of mutually (strongly) embeddable graphs is infinite. On the other hand, there are examples of graphs for which the above number is $1$, for instance a graph that consists of a vertex with finitely many infinite paths originating from it. Thus, a natural question to ask is whether for some graph this number can be finite, but strictly larger than $1$. For strong embeddings this question was first raised by Bonato and Tardif~\cite{tardif}, and is, to our knowledge, still unsolved.

In~\cite{bonato} Bonato and Tardif formulated a similar question for infinite trees. Note that for two trees $T$ and $T'$, an embedding of $T$ into $T'$ is automatically strong, and therefore $T \sim T'$ implies $T \approx T'$. On the other hand, if a tree $T$ and a graph $G$ are mutually embeddable or strongly embeddable, $G$ need not be connected, and so it can happen that $T\sim G$, but not $T \approx G$. In this paper, however, we consider only embeddings between trees, therefore we can safely write ``mutually embeddable'' or ``twinned'' meaning $T \approx T'$. All trees we deal with are infinite, although we shall occasionally emphasise this. We work entirely within ZFC.

Like general graphs, twinned trees need not be isomorphic. For example, let $T$ be a one-way infinite path with a leaf attached to each vertex. Then the removal of any finite number of leaves yields a tree $T'$ with $T'\approx T$, but in general $T$ and $T'$ are non-isomorphic. In fact, if we remove the first $n$ attached leaves for each $n \in \mathbbm{N}$, we obtain infinitely many pairwise non-isomorphic trees $T'$ with $T'\approx T$. 

Given a tree $T$, define the \emph{twin number} of $T$, written $m(T)$, to be the cardinality of the set of isomorphism classes of trees $T'$ with $T'\approx T$. The above example, as well as many others, prompted Bonato and Tardif~\cite{bonato} to make the following \emph{tree alternative conjecture}.

\begin{conjecture} \label{conj:treealt}
For every tree $T$, $m(T)$ is either $1$ or $\infty$.
\end{conjecture}

Note that we do not distinguish here between infinite cardinalities, although even for a countable tree the twin number can be both countable and uncountable. Bonato and Tardif~\cite{bonato} proved their conjecture for rayless trees, i.e.\ for trees without one-way infinite paths. Their main idea was to show an analogous result for rooted rayless trees and then to extend it to general rayless trees by applying a theorem of Halin~\cite{halin}. 

The analogue of Conjecture 1 for rooted trees can be formulated as follows. A rooted tree $(T,r)$ is a tree $T$ with a distinguished vertex $r$, called the \emph{root}. Two rooted trees $(T,r)$ and $(T',r')$ are said to be \emph{twins} or \emph{mutually embeddable}, in notation $(T,r)\approx (T',r')$, if there are injective graph homomorphisms $\phi \colon T \rightarrow T'$ and $\psi \colon T' \rightarrow T$ such that $\phi(r)=r'$ and $\psi(r')=r$. Similarly, $(T,r)$ and $(T',r')$ are isomorphic if there is an isomorphism between $T$ and $T'$ that maps $r$ to $r'$. As above, for a given rooted tree $(T,r)$, let the \emph{twin number} $m(T,r)$ be the number of isomorphism classes of rooted trees $(T',r')$ with $(T',r')\approx (T,r)$. Again, we only distinguish between natural numbers and $\infty$. The authors of~\cite{bonato} implicitly conjectured the following:

\begin{conjecture} \label{conj:rootedtreealt}
Every rooted tree $(T,r)$ has twin number either $1$ or $\infty$.
\end {conjecture}

As was already mentioned, for rayless trees Conjecture~\ref{conj:rootedtreealt} was proved by Bonato and Tardif~\cite{bonato}. Our aim in this paper is to prove it in full.

\begin{theorem} \label{thm:rootedtreealt}
Conjecture~\ref{conj:rootedtreealt} holds.
\end{theorem}

In the next section we establish some basic properties of rooted trees and give a proof of Theorem~\ref{thm:rootedtreealt}. In Section~\ref{sec:locfinitetrees} we discuss Theorem~\ref{thm:rootedtreealt} and Conjecture~\ref{conj:treealt} applied to locally finite trees, i.e.\ trees without vertices of infinite degree. For locally finite trees Theorem~\ref{thm:rootedtreealt} is particularly easy to prove and can be strengthened to $m(T,r)=1$. We also make some progress towards proving Conjecture~\ref{conj:treealt} for locally finite trees. Finally, in Section~\ref{sec:remarks} we formulate two new conjectures on embeddings of locally finite trees and raise several other questions.

\section{Proof of Theorem~\ref{thm:rootedtreealt}}
\label{sec:mainproof}

Let $(T,r)$ be a rooted tree. Denote the neighbours of $r$ by $(r_i)_{i = 1}^{\alpha}$ , where $\alpha$ is an ordinal. Let $T_i$ be the connected component of $T-\left\{rr_i\right\}$ that contains $r_i$. We shall call $(T_i,r_i)$ a \emph{branch} of $(T,r)$. The following lemma from~\cite{bonato} provides a very useful ``recursive'' tool for dealing with rooted trees.

\begin{lemma} \label{lem:bt1a}
If in a rooted tree $(T,r)$ some branch $(T_i,r_i)$ is twinned with infinitely many pairwise non-isomorphic trees, then the same applies to $(T,r)$ itself.
\end{lemma}

Let us slightly extend the statement of Lemma~\ref{lem:bt1a}. The next lemma shows that $m(T,r)$ is at least as large as the maximal $m(T_i,r_i)$ over all branches $(T_i,r_i)$. This corresponds exactly to the statement of Lemma~\ref{lem:bt1a} if $m(T_i,r_i)=\infty$ for some $i$. Moreover, we show that if $m(T,r)<\infty$, then equality between the maximum over all branches and $m(T,r)$ can be attained only if all remaining branches $(T_j,r_j)$ satisfy $m(T_j,r_j)=1$, i.e.\ if none of them can be mutually embedded into any other tree. 

\begin{lemma} \label{lem:bt1b}
For any rooted tree $(T,r)$ we have $m(T,r) \geq \max_{1\leq i \leq \alpha}m(T_i,r_i)$. Furthermore, if for some $i$ we have $1<m(T,r)=m(T_i,r_i)<\infty$, then $m(T_j,r_j)=1$ for all $j\neq i$.
\end{lemma}

\begin{proof}
Suppose that $m(T_i,r_i) \geq n$ for some $i\leq \alpha$ and $n \in \mathbbm{N}$.  Let $(T_i^1,r')$, $(T_i^2,r'), \dotsc, (T_i^n,r')$ be pairwise non-isomorphic trees which are twinned with $(T_i,r_i)$. Let $(T^k,r)$ be the tree obtained from $(T,r)$ by replacing all branches $(T_{i'},r_{i'})\approx (T_i,r_i)$ with a copy of $(T_i^k,r')$. Then the resulting trees $(T^1,r),(T^2,r), \dotsc,(T^n,r)$ are twinned with $(T,r)$, but pairwise non-isomorphic. This implies that $m(T,r)\geq n$, proving the first assertion.

Suppose that in addition $m(T_j,r_j)\geq 2$ for some $j\ne i$. The trees $(T_i,r_i)$ and $(T_j,r_j)$ are either mutually embeddable or not. Let us consider these two cases separately.

\textbf{Case 1.} $(T_j,r_j)\approx (T_i,r_i)$. Replace all branches $(T_{i'},r_{i'})\approx (T_i,r_i)$ except $(T_j,r_j)$ with a copy of $(T_i^1,r')$ and replace $(T_j,r_j)$ with a copy of $(T_i^2,r')$. This gives another tree, which is twinned with $(T,r)$ but not isomorphic to any of the previously constructed.

\textbf{Case 2.} $(T_j,r_j) \not \approx (T_i,r_i)$. Since $m(T_j,r_j)\geq 2$, there must exist a rooted tree $(T'_j,r'_j)$ such that $(T'_j,r'_j)\approx (T_j,r_j)$ and $(T'_j,r'_j)\ncong (T_j,r_j)$. Now take $(T^1,r)$ and replace all branches that are twinned with $(T_j,r_j)$, including $(T_j,r_j)$ itself, by a copy of $(T'_j,r'_j)$. Again, we obtain a new isomorphism class of trees twinned with $T$.

Therefore, we have $m(T,r)>m(T_i,r_i)$ in both cases.
\end{proof}

With Lemma~\ref{lem:bt1b} at our disposal we can find a lower bound on $m(T,r)$ if there is a branch $m(T_i,r_i)$ with $m(T_i,r_i)>1$. However, if for all branches we have $m(T_i,r_i)=1$, Lemma~\ref{lem:bt1b} does not give us any information. The next lemma, which allows us to deal with this case, was also proved by Bonato and Tardif~\cite{bonato}.

\begin{lemma} \label{lem:bt2}
If all branches $(T_i, r_i)$ of a rooted tree $(T,r)$ satisfy $m(T_i, r_i) = 1$, then $m(T,r)$ is $1$ or $\infty$.
\end{lemma}

From here it is a short step to establishing the result of Bonato and Tardif for rooted rayless trees. Indeed, if $T$ is rayless and $1<m(T,r)<\infty$, then by Lemma~\ref{lem:bt1a} and Lemma~\ref{lem:bt2} we must have $1<m(T_i,r_i)<\infty$ for some branch. Now apply the above lemmas to $(T_i,r_i)$ and repeat. Since the tree is rayless, the procedure must stop at some point. That means we would obtain $1<m(T_0,r_0)<\infty$ for a tree $T_0$ consisting of a single vertex $r_0$, clearly a contradiction.

In order to prove Theorem~\ref{thm:rootedtreealt}, we apply Lemma~\ref{lem:bt1b} rather than Lemma~\ref{lem:bt1a}.

\begin{proofthm1}
Suppose that $1 < m(T,r)=n< \infty$. For each vertex $v$ of $T$ define $T(v)$ to be the tree spanned by $v$ and its descendants 
in the tree-order of $(T,r)$. Abusing notation, we write $m(v)$ for $m(T(v),v)$. In particular, $m(r)=m(T,r)$.
By Lemma~\ref{lem:bt1b}, $m(v)\le n$ for each $v$. This means that by passing to an appropriate subtree and perhaps to a different $n>1$, we may assume that for any $v \in T$ the value $m(v)$ is either $1$ or $n$. Under this assumption it follows from Lemma~\ref{lem:bt1b} and Lemma~\ref{lem:bt2} that $T$ contains a ray $P := \left\{r,v_1,v_2,...\right\}$ such that $m(v_i)=n$ for all $i$ and $m(v)=1$ for all $v \notin P$. 

Let $(T'_1,r'),(T'_2,r'),\dotsc,(T'_n,r')$ be $n$ pairwise non-isomorphic trees twinned with $(T(v_1),v_1)$. Let $(T^i,r)$ be the tree obtained from $(T,r)$ by replacing $(T(v_1),v_1)$ with a copy of $(T'_i,r')$. For each $i$, $(T^i,r)$ is twinned with $(T,r)$ by fixing $T-T(v_1)$ and perturbing $(T(v_1),v_1)$. On the other hand, if there is an isomorphism $\phi \colon (T^i,r)\rightarrow (T^j,r)$, then $v_1$ must be mapped on itself, for in both trees $v_1$ is the only neighbour $v$ of $r$, with $m(v)=n$. But then the restriction of $\phi$ to $(T^i(v_1),v_1)$ yields an isomorphism between $(T'_i,r')$ and $(T'_j,r')$, a contradiction. Hence, the trees $(T^i,r)$ represent all $n$ isomorphism classes of trees twinned with $(T,r)$.
It follows that whenever $(T',r')$ satisfies $(T',r')\approx (T,r)$, the root $r'$ must have a neighbour $v'_1$ such that $(T'(v'_1),v'_1)\approx (T(v_1),v_1)$ and $(T'-T'(v'_1),r')\cong (T-T(v_1),r)$.

Since $(T(v_1),v_1)$ also contains a ray $P_1 := P-r = \left\{v_1,v_2,...\right\}$, we can apply the above argument to $(T(v_1),v_1)$ to obtain that whenever $(T',r')\approx (T,r)$, there is a successor of $v'_1$, say $v'_2$, such that $(T'(v'_2),v'_2)\approx (T(v_2),v_2)$ and $(T'-T'(v'_2),r')\cong (T-T(v_2),r)$. Applying the argument repeatedly to all $v_i \in P$ we obtain an isomorphism between $(T',r')$ and $(T,r)$. This implies that $m(T,r)=1$, clearly a contradiction.
\end{proofthm1}

\section{Locally finite trees}
\label{sec:locfinitetrees}

Locally finite trees are in some sense the simplest of all infinite trees. They are countable and, by K\"onig's lemma, a tree which is locally finite and rayless must be finite. Therefore, it seems plausible that a proof of Conjecture~\ref{conj:treealt} for locally finite trees could be a first step towards proving it in full. 
The following lemma shows that for \emph{rooted} locally finite trees there is only one possible value of $m(T,r)$, namely $m(T,r)=1$.

\begin{lemma} \label{lem:loc1}
If $T$ is a locally finite tree, then any embedding, i.e.\ injective homomorphism, of $(T,r)$ into itself is surjective.
\end{lemma}

\begin{proof}
Let $\varphi \colon (T,r)\rightarrow (T,r)$ be an embedding. Note that $\varphi$ preserves the distance of a vertex from $r$, i.e.\ $d(r,v)=d(r,\varphi(v))$. Therefore, for each $n \in \mathbbm{N}$, $\varphi$ induces a self-embedding of the subtree $(T_n,r)$ spanned by $\left\{v \in T : d(r,v)\le n\right\}$. But since the latter trees are finite, $\varphi$ is surjective on each of them, whence, $\varphi$ is surjective on the whole of $(T,r)$.
\end{proof}

\begin{corollary}
If $T$ is a locally finite tree, then $m(T,r)=1$.
\end{corollary}

Indeed, suppose that $\phi \colon (T,r) \rightarrow (T',r')$ and $\psi \colon (T',r') \rightarrow (T,r)$ are injective homomorphisms. Then $\psi \circ \phi \colon (T,r) \rightarrow (T,r)$ is a self-embedding. By Lemma~\ref{lem:loc1} it must be surjective. Thus, $\psi$ is surjective. Since a bijective homomorphism between two trees is an isomorphism, $(T,r)\cong(T',r')$ holds.

Our next aim is to attack Conjecture~\ref{conj:treealt} for locally finite trees. While, unfortunately, we cannot claim to have a proof, we believe that the following observations could be of great help. First we provide a simple but useful isomorphism criterion for rooted locally finite trees. It was first proved by Halin~\cite{halin2} but, for the sake of completeness, we recall the proof. We use the notation from the proof of Lemma~\ref{lem:loc1}.

\begin{lemma} \label{lem:lociso}
Two locally finite rooted trees $(T,r)$ and $(T',r')$ are isomorphic if they are \emph{locally isomorphic} , i.e.\ $(T_n,r)\cong (T'_n,r')$ for all $n \in \mathbbm{N}$.
\end{lemma}

Note that this statement does not generalise even to countable trees --- let $(T,r)$ be the tree consisting of $r$ with countably many finite paths of \emph{each} finite length attached to it and define $T'$ to be $T$ with one additional infinite path attached to $r$. Then $(T,r)\ncong (T',r)$ as the former is rayless and the latter is not. On the other hand $(T_n,r)$ and $(T'_n,r)$ are isomorphic, because both comprise countably many paths of each length $m \leq n$ attached to $r$.

\begin{proof}
This is a standard compactness argument. The number of isomorphisms between two rooted finite trees is always finite. Also, for $m<n$ an isomorphism $ \varphi_n \colon (T_n,r) \rightarrow (T'_n,r')$ induces an isomorphism between $(T_m,r)$ and $(T'_m,r')$ by restriction. Thus, passing to an appropriate subsequence of local isomorphisms, we may assume that the restriction of any $\varphi_n$ on $(T_1,r)$ yields $\varphi_1$. Now we can pass to the next subsequence and assume that for $n\geq 2$ we have $\varphi_n|_{_{T_2}} = \varphi_2$.
Repeating this procedure $n$ times for each $n$, we obtain a sequence of \emph{nested} isomorphisms $\varphi_n \colon (T_n,r) \rightarrow (T'_n,r')$. Finally, since they are nested, we can ``put them together'', i.e.\ define $\varphi(v) = \varphi_n(v)$, where $n = d_{T}(r,v)$, to obtain an isomorphism $\varphi \colon (T,r)\rightarrow (T',r')$.
\end{proof}

Turning to unrooted isomorphisms, we can now formulate the following criterion.

\begin{corollary}
Two locally finite trees $T$ and $T'$ are isomorphic if there exists an assignment of some vertex $r\in T$ to some $r'\in T'$ such that the rooted trees $(T,r)$ and $(T',r')$ are locally isomorphic.
\end{corollary}

Note that once such an assignment exists for some $r \in T$, it exists for any other vertex of $T$ as well. In other words, if we want to show that two locally finite trees $T$ and $T'$ are not isomorphic, it suffices to fix a vertex $r \in T$ and show that for no $r' \in T'$ the rooted trees $(T,r)$ and $(T',r')$ are locally isomorphic. 

Now we are ready to prove Conjecture~\ref{conj:treealt} for a fairly large class of locally finite trees.

\begin{theorem} \label{thm:infcomp}
If $T$ is a locally finite tree, $S\approx T$ and there exists an embedding $\phi \colon S \rightarrow T$ such that $T-\phi(S)$ has infinitely many components, then $m(T)=\infty$.
\end{theorem}

\begin{proof}
Suppose, for a contradiction, that $m(T)=n<\infty$ and let $T^1,T^2,\dotsc,T^n$ be representatives of the isomorphism classes. Note that since $S \approx T$, every tree $T'$ with the ``sandwich property'' $\phi(S)\subset T' \subset T$ is twinned with $T$. So our task is to find a tree having the sandwich property but being not isomorphic to any $T^i$. For this sake we use the above corollary.

Fix a vertex $r \in \phi(S)$ and list all pairs, i.e.\ rooted trees, $(T^i, r')$ where $r' \in T^i$ and ranges over all vertices of $T^i$, for all $i$. Note that this list is countable. Now work through the list, starting with the first pair $(T^{i_1},r_1)$ and check, whether it is isomorphic to $(T,r)$. If it is, adjust $T'=\phi(S)$ by adding a component of $T-\phi(S)$ in the same way as it lies in $T$. Then the local isomorphism breaks down at some $n_1$ and so does the global one. Remember $n_1$. If no, then by Theorem~\ref{thm:infcomp} the local isomorphism already fails at some $n_1$, which we then remember.

Take the next pair $(T^{i_2}, r_2)$ and check, whether $(T^{i_2}, r_2)\cong(T,r)$. If yes, adjust $T'$ by adding another component of $T-\phi(S)$, which is connected to a vertex $v \in \phi(S)$ with $d(r,v)>n_1$ --- this is always possible by the local finiteness and the fact that the number of components is infinite. So the local isomorphism for $(T^{i_1},r_1)$ and $(T,r)$ still fails at $n_1$ and the one for $(T^{i_2},r_2)$ and $(T,r)$ fails now at some $n_2 > n_1$. Remember $n_2$. If no, then the local isomorphism already fails for some $n_2>n_1$, which we then remember.

Repeating this procedure infinitely often we can disturb the local isomorphism for each pair $(T^i,r')$ from the list. Since $T'$ has the sandwich  property at each step of the construction, the limit tree $T'$ is well-defined and has the sandwich property as well. And since we chose the sequence $n_1,n_2,\dotsc$ to be increasing, we must have $(T',r)\ncong (T^i,r')$ for any pair from the list.
\end{proof}

Remarkably, the proof actually shows that $m(T)$ is uncountable, since countably many representatives $T^1,T^2,\dotsc$ also would give rise to a countable list of rooted trees.

\begin{corollary}
If $T$ is a locally finite tree and has a self-embedding $\phi \colon T \rightarrow T$ such that $T-\phi\left[T\right]$ has infinitely many components, then $m(T)=\infty$.
\end{corollary}

So we know that if $T-\phi(S)$ consists of infinitely many components, then $m(T)= \infty$, but what can we say about the components themselves? Using Theorem~\ref{thm:infcomp}, it is not hard to pose a serious restriction on the ``size'' of the components.

Let us define a \emph{nearly finite} tree to be a finite tree with finitely many rays attached to it. 
Equivalently, $T$ is nearly finite if it is locally finite and has only finitely many vertices of degree $3$ or more. The next lemma gives another equivalent characterization of nearly finite trees. 

\begin{lemma} \label{lem:nfinite}
A tree $T$ is nearly finite if and only if it is locally finite and contains no comb (cf.~\cite{diestel}), i.e.\ no ray with infinitely many disjoint paths of finite length (or equivalently, of length~$1$) attached to it, as a subtree.
\end{lemma}

\begin{proof}
The ``only if'' part is obvious. To show the ``if'' part, let us assume that $T$ is not nearly finite and consider the connected hull $S$ of the set of all vertices of degree at least $3$. Either $S$ has a vertex of infinite degree, in which case $T$ is not locally finite, or, by K\"onig's lemma, $S$ contains a ray, which means $T$ contains a comb as a subtree.
\end{proof}

\begin{theorem} \label{thm:nfinite}
If $T$ is a locally finite tree, $S\approx T$ and there exists an embedding $\psi \colon S \rightarrow T$ such that a component of $T-\psi(S)$ is not nearly finite, then $m(T)=\infty$. 
\end{theorem}

\begin{proof}
If a component of $T-\psi(S)$ is not nearly finite, it must by Lemma~\ref{lem:nfinite} contain a comb $C$. Let $T'$ be the connected hull of the spine of $C$ together with $\psi(S)$. Then $\psi(S)\subset T' \subset T$ so $T\approx T'$ but $T'$ has infinitely many components, for it contains the spine of $C$ but at most one of its teeth. Therefore, we are done by Theorem~\ref{thm:infcomp}.
\end{proof}

\begin{corollary}
If $T$ is a locally finite tree and there exists a self-embedding $\psi \colon T \rightarrow T$ such that a component of $T - \psi(T)$ is not nearly finite, then $m(T)=\infty$.
\end{corollary}

Theorem~\ref{thm:nfinite} can also be extended to arbitrary countable trees, which also gives us a new proof that does not rely on Theorem~\ref{thm:infcomp}.

\begin{theorem} \label{thm:countable}
If $T$ is a countable infinite tree, $S\approx T$ and there is an embedding $\psi \colon S \rightarrow T$ such that a component of $T-\psi(S)$ contains a comb, then $m(T)=\infty$.
\end{theorem}

\begin{proof}
Suppose that $m(T)<\infty$ and that a component of $T-\psi(S)$ contains a comb $C$ with all teeth being paths of length $1$. Let $\bar{T}$ be the connected hull of $\psi(S)$ and $C$. Now, altering $\bar{T}$ by removing some teeth of $C$, but leaving infinitely many, we obtain a family $\Theta$ satisfying $\left|\Theta\right|=2^{\aleph_0}$ and $T'\approx T$ for any $T' \in \Theta$. Since $m(T)<\infty$, we know that uncountably many members of $\Theta$ lie in the same isomorphism class. Note however that any tree from $\Theta$ contains an \emph{isolated} comb, i.e.\ a comb with only one exterior edge attached to it, as a subgraph. So some countable $T'\approx T$ contains uncountably many pairwise different isolated combs, all attached to the remaining part of the tree in the same way. But since each vertex of $T'$ can be contained only in countably many isolated combs, we obtain a contradiction to the countability of $T'$.
\end{proof}

\section{Concluding remarks} 
\label{sec:remarks}

We have shown that $m(T)=\infty$ for a fairly large class of locally finite trees $T$. However, we have not managed to find a locally finite tree with $m(T)= 1$ in addition to the trees all whose self-embeddings are surjective, like regular (or $d$-ary) trees, and the ray, which can be considered an exceptional very ``small'' tree. This prompts us to conjecture the following.

\begin{conjecture} \label{conj:loc1}
If $T$ is a locally finite tree and has a non-surjective self-embedding, then $m(T)=\infty$ unless $T$ is a ray.
\end{conjecture}

Obviously, Conjecture~\ref{conj:loc1} implies Conjecture~\ref{conj:treealt} for locally finite trees. Also, one could relax Conjecture~\ref{conj:loc1} by replacing $m(T)=\infty$ with $m(T)>1$.

\begin{conjecture} \label{conj:loc2}
If $T$ is a locally finite tree and has a non-surjective self-embedding, then $m(T)>1$ unless $T$ is a ray.
\end{conjecture}

We believe, however, that this version should be as hard as Conjecture~\ref{conj:loc1}, in other words once we have proved Conjecture~\ref{conj:loc2}, we know all about possible values $m(T)$ and then we either have proved Conjecture~\ref{conj:loc1} or disproved Conjecture~\ref{conj:treealt}, which we think would be even more interesting.

A way to attack Conjecture~\ref{conj:loc2} would be to show that there is no locally finite tree with infinitely many leaves $x_1,x_2, \dotsc$ such that $T \cong T-x_1\cong T-x_2\cong \dotsb$. This would imply Conjecture~\ref{conj:loc2}, but unfortunately the claim is false. In fact, in~\cite{tyomkyn} we show much more --- there exists a tree $T$ of maximal degree $3$ and infinitely many leaves satisfying $T \cong T-x$ for \emph{any} leaf $x$.

Turning to the original definition of a mutual embedding of graphs, let us remark that the number of pairwise non-isomorphic \emph{graphs} twinned with a given tree $T$ is trivially $1$ or $\infty$, regardless of whether $T$ is locally finite or not. Indeed, if $T$ allows for a non-surjective self-embedding, then $T$ is twinned with any graph consisting of a copy of $T$ and a finite number of isolated vertices. The analogous question with ``embeddable'' replaced by ``strongly embeddable'' is harder and seems to be closely related to Conjecture 1. A rather exotic yet appealing question would be to determine possible cardinalities of ``$\approx$''-equivalence classes of all graphs $H$ with $H \sim G$. 

One could try to extend the definition of $m(T)$ to a general connected graph $G$ by considering the number of equivalence classes of all \emph{connected} graphs $H$ with $G \sim H$ or $G \approx H$. Again it would be very interesting to find a graph $G$ for which one of these numbers is finite, but larger than $1$ or to prove that no such graph can exist.

Similar problems arise for other relations of graphs like the minor or the topological minor relation. Recently, Matthiesen~\cite{matthiesen}
proved a result of the above kind for the topological minor relation of trees.

It is likely that these problems, especially the ones about trees, could be successfully tackled by using group theoretical techniques. We can also imagine that some meta-theorems from logic, in particularly model theory, could find an application in this area. 

\section*{Acknowledgements}
\label{sec:acknowledgements}

I would like to thank my PhD supervisor Professor Béla Bollobás for suggesting this problem and for his help and advice.


\begin{thebibliography}{}
\bibitem[1]{tardif} A.Bonato, C.Tardif, Large families of mutually embeddable vertex-transitive graphs, \textit{J. Graph Theory} \textbf{43} (2003) pp. 99--106. 
\bibitem[2]{bonato} A.Bonato, C.Tardif, Mutually embeddable graphs and the tree alternative conjecture, \textit{J. Combinatorial Theory, Series B} \textbf{96} (2006), pp. 874--880.
\bibitem[3]{diestel} R.Diestel, Graph theory 3rd Ed., Springer Verlag (2005).
\bibitem[4]{halin}   R. Halin, Fixed configurations in graphs with small number of disjoint rays, in: R. Bodendiek, Editor, \textit{Contemporary Methods in Graph Theory}, Bibliographisches Inst., Mannheim (1990), pp. 639--649.
\bibitem[5]{halin2} R.Halin, Automorphisms and endomorphisms of infinite locally finite graphs, \textit{Abh. Math. Sem. Univ. Hamburg} \textbf{39} (1973) pp. 251--283.
\bibitem[6]{matthiesen} L. Matthiesen, There are uncountably many topological types of locally finite trees, \textit{J. Combinatorial Theory, Series B} \textbf{96} (2006), pp. 758--760.
\bibitem[7]{rado} R.Rado, Universal graphs and universal functions, \textit{Acta Arithmetica} \textbf{9} (1964) pp. 331--340.
\bibitem[8]{tyomkyn} M.Tyomkyn, A locally finite tree that behaves like an infinite star. Preprint.
\end{thebibliography}
\end{document}